\documentclass[10pt,a4paper]{article}
\usepackage{fullpage}
\usepackage[latin1]{inputenc}
\usepackage[english]{babel}
\usepackage{amsfonts}
\usepackage{amssymb}
\usepackage{amsmath}
\usepackage{latexsym}
\newtheorem{theoreme}{Theorem}
\newtheorem{lemme}[theoreme]{Lemma}
\newtheorem{corollaire}[theoreme]{Corollary}
\newtheorem{proposition}[theoreme]{Proposition}
\newtheorem{defi}[theoreme]{Definition}
\newenvironment{preuve}{\emph{Proof} : }{\begin{flushright}$\Box$\end{flushright}\vspace{0.5cm}}

\newenvironment{preuvep}{\emph{Proof of proposition} : }{\begin{flushright}$\Box$\end{flushright}\vspace{0.5cm}}
\newenvironment{remarque}{\textbf{\emph{Remark}} : }{\vspace{0.5cm}}
\newcounter{numalgo}
\newenvironment{algo}[1]{\stepcounter{numalgo}\vspace{0.2cm}\hrule\vspace{0.1cm}\noindent\textbf{Algorithm  \arabic{numalgo}} #1\vspace{0.1cm}\hrule\vspace{0.1cm}}{\vspace{0.1cm}\hrule\vspace{0.5cm}}
\newcommand{\F}[1]{\mathbb{F}_{#1}}
\title{On the covering radius of first order generalized Reed-Muller codes}
\author{Elodie Leducq\footnote{IMJ, preparing a PhD thesis under direction of Jean-françois Mestre}}
\date{}
\begin{document}

\maketitle

\section{Introduction}
The determination of the covering radius of the first order Reed-Muller code is a difficult problem in coding theory. In this paper, we generalize to any $q$ some results proved in \cite{BER,MYK,HOU,HOU2} for $q=2$.
\\\\Let $q=p^t$, $p$ a prime number.\\ Let $B_m^q=\F{q}[X_1,\ldots,X_m]/(X_1^q-X_1,\ldots,X_m^q-X_m)$; $B_m^q$ actually consists of all the functions from $\F{q}^m$ to $\F{q}$. We identify $B_m^q$ with $\F{q}^{q^m}$ through the application
$$\begin{array}{rcl}B_m^q&\rightarrow&\F{q}^{q^m}\\P&\mapsto & \left(P(x)\right)_{x\in\F{q}^m}\end{array}$$
For all $b\in\F{q}^m$, we denote by $\mathrm{1}_b$ the function in $B_m^q$ such that $\mathrm{1}_b(b)=1$ and for all $x\neq b$, $\mathrm{1}_b(x)=0$.
\\\\The weight $|P|$ of $P\in B_m^q$ is $\mathrm{Card}(\{x, P(x)\neq0\})$. The Hamming distance in $B_m^q$ is denoted by $d(.,.)$.
\\For $0\leq r\leq m(q-1)$, the $r$th order generalized Reed-Muller code of length $q^m$ is $$R_q(r,m)=\{P\in B_m^q, \mathrm{deg}(P)\leq r\}$$ where $\mathrm{deg}(P)$ is the degree of the representant of $P$ with degree at most $q-1$ in each variable (see \cite{HAB}).
\\For all $0\leq r\leq m(q-1)$, the affine group $\mathrm{GA}_m(\F{q})$ acts on $R_q(r,m)$ by its natural action and we have
\begin{proposition}(see \cite{HAB})\\For all $q$, for all $m\geq0$ and $1\leq r\leq m(q-1)-2$, $$Aut(R_q(r,m))=\mathrm{GA}_m(\F{q}).$$\end{proposition}
The covering radius of a code $C$ of length $n$ is $$\rho(C)=\displaystyle\max_{x\in\F{q}^n}\min_{c\in C}|x-c|$$
We denote by $\rho(r,m)$ the covering radius of $R_q(r,m)$.
\\\\For $q=2$, $\rho(1,m)$ are unknown for $m$ odd, $m\geq 9$. We know that $\rho(1,5)=12$ and $\rho(1,7)=56$ (see \cite{BER,HOU}).
\\Since $R_q(1,m)\subset R_q(2,m)\subset\ldots\subset R_q(m(q-1),m)=\F{q}^{q^m}$, we can try to study $\rho(1,m)$ through the covering radius of $R_q(1,m)$ in $R_q(r,m)$, $2\leq r\leq m(q-1)$. We define $$\rho_r(1,m)=\displaystyle\max_{x\in R_q(r,m)}\min_{c\in R_q(1,m)}|x-c|$$
For $q=2$, it is known that $\rho_2(1,m)=2^{m-1}-2^{\lceil \frac{m}{2}\rceil-1}$ (see \cite{HOU2}) which gives $\rho(1,m)\geq 2^{m-1}-2^{\lceil \frac{m}{2}\rceil-1}$
and, since $\rho(1,m)\leq 2^{m-1}-2^{\frac{m}{2}-1}$ (see \cite{MYK}), we get, for $m$ even $\rho(1,m)=2^{m-1}-2^{ \frac{m}{2}-1}$.
\\\\In part \ref{2}, we give a general upper bound for covering radius of codes over $\F{q}$.
\\Then we calculate $\rho_2(1,m)$ :
\begin{theoreme}\label{rho2}For all $q$, $m\geq0$, $$\rho_2(1,m)=(q-1)q^{m-1}-q^{\lceil\frac{m}{2}\rceil-1}.$$\end{theoreme}
We use previous results to give an upper bound and a lower bound for $\rho(1,m)$, which give $\rho(1,m)$ when $m$ is even :
\begin{theoreme}\label{bound}For all $q$, $m\geq0$,
$$(q-1)q^{m-1}-q^{\lceil\frac{m}{2}\rceil-1}\leq\rho(1,m)\leq (q-1)q^{m-1}-q^{\frac{m}{2}-1}.$$ \end{theoreme}
Furthermore, we get functions, $f$, such that $d(f,R(1,m))=\rho(1,m)$ when $m$ is even.
\\Finally, we study more precisely the case where $q=3$.

\section{A general upper bound}\label{2}
We need to extend the definition of self-complementary code and the definition of strength of a code for a binary code (see \cite{MYK}) to a code over $\F{q}$.

\begin{defi}A code $C$ over $\mathbb{F}_q$ is self-complementary if $$\forall c\in C,\ \forall \omega\in\mathbb{F}_q, \ \overline{c}^{\omega}=(\omega,\ldots,\omega)+c\in C.$$ \end{defi}

\begin{defi} A code $C$ over $\mathbb{F}_q$ has strength $s$ if each $s$-subset of coordinates of the code contains all elements of $\mathbb{F}_q^s$ a constant number of times.\end{defi}

Now, we generalize the upper bound of covering radius of a binary code given in \cite{MYK} to codes over $\F{q}$.

\begin{lemme} \label{force}Let $C$ be a code over $\mathbb{F}_q$ of length $n$, and let $v\in\mathbb{F}_q^m$.
\\If $C$ has strength 2, then
$$\sum_{u\in v+C}|u|^2=n\left(\frac{q-1}{q}\right)\left((n-1)\left(\frac{q-1}{q}\right)+1\right)\mathrm{Card}(C)$$\end{lemme}

\begin{preuve} $v+C$ has strength 2 if and only if $C$ has strength 2.
\\Furthermore, $\mathrm{Card}(v+C)=\mathrm{Card}(C)$, so it is enough to prove the lemma for $v=0$.
\begin{eqnarray*}\sum_{u\in C}|u|^2&=&\sum_{u\in C}\sum_{i, u_i\neq 0}\sum_{j,u_j\neq 0} 1 = \sum_{u\in C}\sum_{i,u_i\neq0}\sum_{j\neq i,u_j\neq0} 1+\sum_{u\in C}\sum_{i,u_i\neq0} 1\\&=&\sum_{i=1}^n\sum_{j\neq i}\sum_{u\in C,u_iu_j\neq0}1+\sum_{i=1}^n\sum_{u\in C,u_i\neq0}1\end{eqnarray*}
Since $C$ has strength 2, $\displaystyle\sum_{u\in C,u_iu_j\neq0}1=\left(\frac{q-1}{q}\right)^2\mathrm{Card}(C)$
\\and $\displaystyle\sum_{u\in C,u_i\neq0}1=\left(\frac{q-1}{q}\right)\mathrm{Card}(C)$.
\\Hence $\displaystyle\sum_{u\in C}|u|^2=n(n-1)\left(\frac{q-1}{q}\right)^2 \mathrm{Card}(C)+n\left(\frac{q-1}{q}\right) \mathrm{Card}(C)$. \end{preuve}

\begin{theoreme}\label{gen}If $C$ is a self-complementary code over $\F{q}$ of length $n$ and strength 2, then
$$\rho(C)\leq \frac{(q-1)}{q}n-\frac{\sqrt{n}}{q}.$$\end{theoreme}

\begin{preuve}Let $v\in\F{q}^n$ such that its distance to any codeword is at least $r$, i.e. $\forall u\in v+C$,\ $|u|\geq r$.
\\Since $C$ is self-complementary, if $u\in v+C$ then $\overline{u}^{\omega}\in v+C$ and $|\overline{u}^{\omega}|\geq r$.
\\We have :
$$\sum_{\omega\in \mathbb{F}_q}|\overline{u}^{\omega}|=\sum_{\omega\in\mathbb{F}_q}(|u|-x_{-\omega}+x_0)=q(|u|+x_0)-\sum_{\omega\in\mathbb{F}_q}x_{-\omega}=(q-1)n$$ where
$x_{\omega}=\mathrm{Card}(\{i,u_i=\omega\})$.
\\\\Assume $r>\frac{q-1}{q}n$, so
$$\frac{q-1}{q}n<r\leq|u|=n(q-1)-\sum_{\omega\in\mathbb{F}_q^*}|\overline{u}^{\omega}|\leq n(q-1)-(q-1)r<\frac{q-1}{q}n.$$
We get a contradiction. So we write $r=\frac{q-1}{q}n-\rho$ with $\rho\geq 0$ and we have :
\begin{eqnarray*}\sum_{\omega\in\mathbb{F}_q}|\overline{u}^{\omega}|^2&=&\sum_{\omega\in\mathbb{F}_q^*}|\overline{u}^{\omega}|^2+|u|^2 \\&=&\sum_{\omega\in\mathbb{F}_q^*}|\overline{u}^{\omega}|^2+\left((q-1)n-\sum_{\omega\in\mathbb{F}_q^*}|\overline{u}^{\omega}|\right)^2 \\&\leq&(q-1)r^2+(q-1)^2(n-r)^2 \\&=&(q-1)\left(\frac{q-1}{q}n-\rho\right)^2+(q-1)^2\left(\frac{n}{q}+\rho\right)^2 \\&=&q\left(n^2\frac{(q-1)^2}{q^2}+(q-1)\rho^2\right)\end{eqnarray*}
we get $$\displaystyle\sum_{u\in v+C}|u|^2\leq\left(n^2\frac{(q-1)^2}{q^2}+(q-1)\rho^2\right)\mathrm{Card}(C).$$
And, from lemma \ref{force}, $$n^2\frac{(q-1)^2}{q^2}+n\frac{(q-1)}{q^2}\leq\left(n^2\frac{(q-1)^2}{q^2}+(q-1)\rho^2\right).$$
Hence $\frac{n}{q^2}\leq\rho^2$, and so, since $\rho\geq0$, $r\leq\frac{(q-1)}{q}n-\frac{\sqrt{n}}{q}$.
\\Theorem \ref{gen} follows from the definition of covering radius.
 \end{preuve}

\section{Counting zeros of quadratic forms}

\begin{defi} We say that an application from $\F{q}^m$ to $\F{q}$ is a quadratic form if
\begin{enumerate}\item $Q(ax)=a^2Q(x)$ for $a\in\F{q}$ and $x\in\F{q}^m$.
\item $\phi : \F{q}^m\times\F{q}^m\rightarrow\F{q}$ such that $\phi(x,y)=Q(x+y)-Q(x)-Q(y)$ is a bilinear form.\end{enumerate}
$\phi$ is called the bilinear form associated to $Q$.\end{defi}

\begin{defi}We called rank of a bilinear form $\phi$, the rank of the following application
$$\psi_{\phi} : \begin{array}{rcl}\F{q}^n&\rightarrow&\F{q}^{m*}\\x&\mapsto&(y\mapsto\phi(x,y))\end{array}$$ Let $Q$ be a quadratic form over $\F{q}^n$ and $\phi$ the associated bilinear form.
\\Let $V=\{x\in\mathrm{Ker}(\psi_{\phi}), Q(x)=0\}$ and let $v=\mathrm{dim}(V)\leq n-\mathrm{Rg}(\phi)$, then we define the rank of $Q$ by $\mathrm{Rg}(Q)=n-v$.
\\We say that $Q$ is non degenerate, if $v=0$.\end{defi}

\begin{remarque}If $q$ is odd, then $v=n-\mathrm{Rg}(\phi)$ and $\mathrm{Rg}(\phi)=\mathrm{Rg}(Q)$.\end{remarque}

We need the following theorem about reduction of quadratic forms :

\begin{theoreme}\label{red}Let $Q$ be a quadratic form of rank $R$ on $\F{q}^n$, then there exists a basis $(e_i)_{1\leq i\leq n}$ of $\F{q}^n$ such that :
\begin{itemize}\item If $R=2s+1$, then $Q(\displaystyle\sum_{i=1}^nx_ie_i)=\displaystyle\sum_{i=1}^sx_{2i-1}x_{2i}+ax_{2s+1}^2$ (1)
\item If $R=2s$, then $Q(\displaystyle\sum_{i=1}^nx_ie_i)=\displaystyle\sum_{i=1}^sx_{2i-1}x_{2i}$ (2) \\  or $Q(\displaystyle\sum_{i=1}^nx_ie_i)=\displaystyle\sum_{i=1}^sx_{2i-1}x_{2i}+ax_{2s-1}^2+bx_{2s}^2+cx_{2s-1}x_{2s}$ where $ax^2+cx+b$ is irreducible over $\F{q}$ (3)\end{itemize}\end{theoreme}

\begin{preuve}Let $\phi$ be the bilinear form associated to $Q$.
\\Let $N=\{x, \forall y\in \F{q}^n, \phi(x,y)=0\}$ and $S$ be a subspace supplementary to $N$ in $\F{q}^n$.
\\Let $V_0=\{x\in N,Q(x)=0\}$ and let $V_1$ be a subspace supplementary to $V_0$ in $N$.
\\So $Q$ restricted to $V_1\oplus S$ is non degenerate.
\\\\If $q$ is even, see \cite{DIE} p.33-34 and \cite{dic} p.197-199.
\\\\If $q$ is odd, see \cite{HIR} p.117-118 and p.121-123.\end{preuve}

Now, we can count zeros of quadratic forms (see \cite{HIR}) :

\begin{theoreme}\label{zero}Let $Q$ be a quadratic form over $\mathbb{F}_q^n$ of rank $R$, then the number of zeros of $Q$ is $$N(Q)=q^{n-1}+(\omega-1)(q-1)q^{n-\frac{R}{2}-1}$$
where $\omega=\left\{\begin{array}{ll}1&\textrm{if $R$ is odd}\\2&\textrm{if $R$ is even and $Q$ is of type $(2)$ (see theorem \ref{red})} \\0&\textrm{if $R$ is even and  $Q$ is of type $(3)$}\end{array}\right.$
\end{theoreme}

\begin{preuve}
If $R=0$, there are $q^n$ zeros.
\\\\If $R=1$, we can write $Q=ax_1^2$, and so $Q$ has $q^{n-1}$ zeros.
\\\\If $R=2$, we can write $Q=ax_1^2+bx_2^2+cx_1x_2$.
\\If $ax^2+cx+b$ is irreducible, we are in case $(3)$ and $Q$ has $q^{n-2}=q^{n-1}-(q-1)q^{n-\frac{2}{2}-1}$ zeros. Otherwise, $Q$ factors, so we are in case $(2)$, and $Q$ has $(2q-1)q^{n-2}=q^{n-1}+(q-1)q^{n-\frac{2}{2}-1}$ zeros.
\\\\If $R\geq 3$, by theorem \ref{red}, we can write $Q=x_1x_2+Q^{(1)}(x_3,\ldots,x_n)$ where $Q^{(1)}$ is a quadratic form of rank $R-2$.
\\If $Q^{(1)}(a_3,\ldots,a_n)=0$, then there are $(2q-1)$ couples $(x_1,x_2)$ such that $Q(x_1,x_2,a_3,\ldots,a_n)=0$.
\\Otherwise, there are $(q-1)$ couples $(x_1,x_2)$ such that $Q(x_1,x_2,a_3,\ldots,a_n)=0$.
\\Hence $$N(Q)=(2q-1)N(Q^{(1)})+(q-1)(q^{n-2}-N(Q^{(1)}))=qN(Q^{(1)})+(q-1)q^{n-2}.$$
Continuing this process, we get, for $r$ such that $R-2r\geq 1$
\begin{eqnarray*}N(Q)&=&q^rN(Q^{(r)})+(q-1)(q^{n-2}+q^{n-3}+\ldots+q^{n-(r+1)})\\&=&q^rN(Q^{(r)})+q^{n-1}-q^{n-(r+1)}\end{eqnarray*}
where $Q^{(r)}$ is a quadratic form in $x_{2r+1},\ldots,x_n$ with rank $R-2r$.
\\If $R$ is odd, we put $R=2s+1$ and $r=s$, then we get $N(Q)=q^sq^{n-2s-1}+q^{n-1}-q^{n-s-1}=q^{n-1}$ since $Q^{(s)}$ is a quadratic form in $n-2s$ variables, with rank 1. That gives the theorem in the case where $R$ is odd.
\\If $R$ is even, we put $R=2s$ and $r=s-1$, then 
$Q^{(s-1)}$ is a quadratic form in $n-2s+2$ variables of rank 2.
If $Q$ is of type $(3)$, $Q^{(s-1)}$ does not factor and has $q^{n-2s}$ zeros.
So \\$N(Q)=q^{s-1}q^{n-2s}+q^{n-1}-q^{n-s}=q^{n-1}-(q-1)q^{n-s-1}$.
If $Q$ is of type $(2)$, $Q^{(s-1)}$ factors and has $(2q-1)q^{n-2s}$ zeros.
So $N(Q)=q^{s-1}(2q-1)q^{n-2s}+q^{n-1}-q^{n-s}=q^{n-1}+(q-1)q^{n-s-1}$.
\end{preuve}

Now we are able to prove theorem \ref{rho2}.

\section{Proof of theorem \ref{rho2}}

Let $q_0=\displaystyle\sum_{1\leq i\leq j\leq m}a_{i,j}x_ix_j$
\\In order to get the weight of $q_0+\alpha_1x_1+\ldots+\alpha_mx_m+\beta$, we homogenize : \\let $Q=q_0+\alpha_1x_1z+\ldots+\alpha_mx_mz+\beta z^2$.
\\We denote by $N_q^{\infty}$, the number of zeros of $q_0$ in $\F{q}^m$, which is the number of infinite points of the quadric defined by $Q=0$ and by $N_q$, the number of zeros of $Q$ in $\F{q}^{m+1}$.
\\Then the number of zeros of $q_0+\alpha_1x_1+\ldots+\alpha_mx_m+\beta$, $N$, is the number of point of the quadric which are not infinite points, so we get $$N=\frac{N_q-N_q^{\infty}}{q-1}.$$
\\By theorem \ref{zero} we have :
$$ N_q^{\infty}=q^{m-1}+q^{m-\frac{r}{2}-1}(q-1)(\omega_{q_0}-1)$$ where $r=\mathrm{rg}(q_0)$ and \\ $\omega_{q_0}=\left\{\begin{array}{ll}1&\textrm{if $r$ is odd}\\2&\textrm{if $r$ is even and $q_0$ is of type $(2)$}\\0&\textrm{if $r$ is even and $q_0$ is of type $(3)$}\end{array}\right.$
\\and $$N_q=q^{m}+q^{m-\frac{R}{2}}(q-1)(\omega_Q-1)$$ where $R=\mathrm{rg}(Q)$ and\\ $\omega_Q=\left\{\begin{array}{ll}1&\textrm{if $R$ is odd}\\2&\textrm{if $R$ is even and $Q$ is of type $(2)$} \\0&\textrm{if $R$ is even and $Q$ is of type $(3)$}\end{array}\right.$
\\\\So, the number of zeros of $q_0+\alpha_1x_1+\ldots+\alpha_mx_m+\beta$ is $$N=q^{m-1}-(\omega_{q_0}-1) q^{m-\frac{r}{2}-1}+(\omega_Q-1)q^{m-\frac{R}{2}}.$$
Hence $$|q_0+\alpha_1x_1+\ldots+\alpha_mx_m+\beta|=(q-1)q^{m-1}+(\omega_{q_0}-1)q^{m-\frac{r}{2}-1}-(\omega_Q-1)q^{m-\frac{R}{2}}.$$
\\\\Then we want to calculate $d(q_0,R_q(1,m))$ : \\
\begin{itemize}\item If $r$ is odd, $\omega_{q_0}=1$, $|q_0+\alpha_1x_1+\ldots+\alpha_mx_m+\beta|=(q-1)q^{m-1}-(\omega_Q-1)q^{m-\frac{R}{2}}$ and
$q_0$ can be reduced to $x_1x_2+\ldots+x_{r-2}x_{r-1}+ax_r^2$ by a linear transformation, which does not change the weight, so we can assume that
\begin{eqnarray*}Q&=&x_1x_2+\ldots+x_{r-2}x_{r-1}+ax_r^2+\alpha_1x_1z+\ldots+\alpha_mx_mz+\beta z^2 \\&=&(x_1+\alpha_2z)(x_2+\alpha_1z)+\ldots+(x_{r-2}+\alpha_{r-1}z)(x_{r-1}+\alpha_{r-2}z)+ax_r^2\\&&\hspace{1,5cm}+z(\alpha_rx_r+\ldots+\alpha_mx_m)+\underbrace{(\beta-\alpha_1\alpha_2-\ldots-\alpha_{r-2}\alpha_{r-1})}_{\theta}z^2
\\&=&(x_1+\alpha_2z)(x_2+\alpha_1z)+\ldots+(x_{r-2}+\alpha_{r-1}z)(x_{r-1}+\alpha_{r-2}z)+ax_r^2\\&&\hspace{6cm}+z(\alpha_rx_r+\ldots+\alpha_mx_m+\theta z)\end{eqnarray*}
If there exists $i>r$ such that $\alpha_i\neq0$, then $R=r+2$ and $\omega_Q=1$.
\\If for all $i>r$, $\alpha_i=0$ :
\\ If $\theta\neq0$, then $R=r+1$ and $\omega_Q=\left\{\begin{array}{ll} 0&\textrm{if $ax^2+\alpha_rx+\theta$ is irreductible}\\2&\textrm{otherwise}\end{array}\right.$
\\ If $\theta=0$, then, if $\alpha_r=0$, $R=r$ and $\omega_Q=1$. Otherwise, $R=r+1$ and $\omega_Q=2$.
\\\\Hence $d(q_0,R_q(1,m))=(q-1)q^{m-1}-q^{m-\frac{r+1}{2}}$.\\
\item If $r$ is even and $\omega_{q_0}=2$, $|q_0+\alpha_1x_1+\ldots+\alpha_mx_m+\beta|=(q-1)q^{m-1}+q^{m-\frac{r}{2}-1}-(\omega_Q-1)q^{m-\frac{R}{2}}$ and $q_0$ can be reduced to $x_1x_2+\ldots+x_{r-1}x_{r}$ by a linear transformation.
\begin{eqnarray*}Q&=&x_1x_2+\ldots+x_{r-1}x_{r}+\alpha_1x_1z+\ldots+\alpha_mx_mz+\beta z^2 \\&=&(x_1+\alpha_2z)(x_2+\alpha_1z)+\ldots+(x_{r-1}+\alpha_{r}z)(x_{r}+\alpha_{r-1}z)\\&&\hspace{0.5cm}+z(\alpha_{r+1}x_{r+1}+\ldots+\alpha_mx_m)+\underbrace{(\beta-\alpha_1\alpha_2-\ldots-\alpha_{r-1}\alpha_{r})}_{\theta}z^2
\end{eqnarray*}
If there exists $i>r$ such that $\alpha_i\neq 0$, $R=R+2$ and $\omega_Q=2$.
\\If for all $i>r$, $\alpha_i=0$ and $\theta=0$, $R=r$ et $\omega_Q=2$.
\\If for all $i>r$, $\alpha_i=0$ and $\theta\neq 0$, $R=r+1$ et $\omega_Q=1$.
\\\\Hence $d(q_0,R_q(1,m))=(q-1)q^{m-1}+q^{m-\frac{r}{2}-1}-q^{m-\frac{r}{2}}$.\\
\item If $r$ is even and $\omega_{q_0}=0$, $|q_0+\alpha_1x_1+\ldots+\alpha_mx_m+\beta|=(q-1)q^{m-1}-q^{m-\frac{r}{2}-1}-(\omega_Q-1)q^{m-\frac{R}{2}}$. By a linear transformation, $q_0$ can be reduced to $x_1x_2+\ldots+x_{r-3}x_{r-2}+ax_{r-1}^2+bx_r^2+cx_{r-1}x_r$ with $ax^2+cx+b$ irreducible.
{\setlength\arraycolsep{1pt}\begin{eqnarray*}Q&=&x_1x_2+\ldots+x_{r-3}x_{r-2}+ax_{r-1}^2+bx_r^2+cx_{r-1}x_r+\alpha_1x_1z+\ldots+\alpha_mx_mz+\beta z^2 \\&=&(x_1+\alpha_2z)(x_2+\alpha_1z)+\ldots+(x_{r-3}+\alpha_{r-2}z)(x_{r-2}+\alpha_{r-3}z)+ax_{r-1}^2\\&&\hspace{1cm}+bx_r^2+cx_{r-1}x_r+z(\alpha_{r-1}x_{r-1}+\ldots+\alpha_mx_m)+\underbrace{(\beta-\alpha_1\alpha_2-\ldots-\alpha_{r-3}\alpha_{r-2})}_{\theta}z^2
\end{eqnarray*}}
If there exists $i>r$ such that $\alpha_i\neq0$, $R=r+2$ and $\omega_Q=0$.
\\Assume that for all $i>r$, $\alpha_i=0$.
\\First, we study the case where $q$ is odd.
\\Since $ax^2+cx+b$ is irreducible, we have $a\neq 0$ and $\Delta=b-\frac{c^2}{4a}\neq0$.
{\setlength\arraycolsep{2pt}\begin{eqnarray*}Q&=&(x_1+\alpha_2z)(x_2+\alpha_1z)+\ldots+(x_{r-3}+\alpha_{r-2}z)(x_{r-2}+\alpha_{r-3}z) \\&&\hspace{0.5cm}+a(x_{r-1}+\frac{c}{2a}x_r+\frac{\alpha_{r-1}}{2a}z)^2+\Delta x_r^2+(\alpha_r-\frac{c\alpha_{r-1}}{2a})x_rz+(\theta-\frac{\alpha_{r-1}^2}{4a})z^2\\&=&(x_1+\alpha_2z)(x_2+\alpha_1z)+\ldots+(x_{r-3}+\alpha_{r-2}z)(x_{r-2}+\alpha_{r-3}z) \\&&\hspace{3cm}+a(x_{r-1}+\frac{c}{2a}x_r+\frac{\alpha_{r-1}}{2a}z)^2+\Delta(x_r+\frac{2a\alpha_r-c\alpha_{r-1}}{4a\Delta}z)^2\\&&\hspace{5cm}+(\theta-\frac{\alpha_{r-1}^2}{4a}-\frac{(2a\alpha_r-c\alpha_{r-1})^2}{16a^2\Delta})z^2\end{eqnarray*}}
If $\theta\neq \frac{\alpha_{r-1}^2}{4a}+\frac{(2a\alpha_r-c\alpha_{r-1})^2}{16a^2\Delta}$, $R=r+1$ and $\omega_Q=1$.
\\If $\theta= \frac{\alpha_{r-1}^2}{4a}+\frac{(2a\alpha_r-c\alpha_{r-1})^2}{16a^2\Delta}$, $R=r$ and $\omega_Q=0$, since $ax^2+cx+b$ is irreducible.
\\\\Then we study the case where $q$ is even.
\\Since $ax^2+cx+b$ is irreducible, $c\neq0$. So we have
{\setlength\arraycolsep{2pt}\begin{eqnarray*}Q&=&(x_1+\alpha_2z)(x_2+\alpha_1z)+\ldots+(x_{r-3}+\alpha_{r-2}z)(x_{r-2}+\alpha_{r-3}z)\\&&+c(x_{r-1}+\frac{\alpha_r}{c}z)(x_r+\frac{\alpha_{r-1}}{c}z)
+\left(\sqrt{a}x_{r-1}+\sqrt{b}x_r+\sqrt{\theta-\frac{\alpha_{r-1}\alpha_r}{c}}z\right)^2\end{eqnarray*}}
If $c^2\theta\neq a\alpha_r^2+b\alpha_{r-1}^2+c\alpha_{r-1}\alpha_r$, $R=r+1$ and $\omega_Q=1$.
\\If $c^2\theta= a\alpha_r^2+b\alpha_{r-1}^2+c\alpha_{r-1}\alpha_r$,
{\setlength\arraycolsep{2pt}\begin{eqnarray*}Q&=&(x_1+\alpha_2z)(x_2+\alpha_1z)+\ldots+(x_{r-3}+\alpha_{r-2}z)(x_{r-2}+\alpha_{r-3}z)\\&&+c(x_{r-1}+\frac{\alpha_r}{c}z)(x_r+\frac{\alpha_{r-1}}{c}z) +\left(\sqrt{a}x_{r-1}+\sqrt{b}x_r+(\frac{\sqrt{a}}{c}\alpha_r+\frac{\sqrt{b}}{c}\alpha_{r-1})z\right)^2
\\&=&(x_1+\alpha_2z)(x_2+\alpha_1z)+\ldots+(x_{r-3}+\alpha_{r-2}z)(x_{r-2}+\alpha_{r-3}z) \\&&\hspace{1cm}+c(x_{r-1}+\frac{\alpha_r}{c}z)(x_r+\frac{\alpha_{r-1}}{c}z)+a(x_{r-1}+\frac{\alpha_r}{c}z)^2+b(x_{r}+\frac{\alpha_{r-1}}{c}z)^2\end{eqnarray*}}
so $R=r$ et $\omega_Q=0$.
\\\\Hence $d(q_0,R_q(1,m))=(q-1)q^{m-1}-q^{m-\frac{r}{2}-1}$.
\end{itemize}

\section{Bounds of $\rho(1,m)$}

We use the general upper bound to find an upper bound to $\rho(1,m)$ :

\begin{proposition}\label{max}For all $q$, $m\geq 0$, we have $$\rho(1,m)\leq (q-1)q^{m-1}-q^{\frac{m}{2}-1}.$$\end{proposition}

\begin{preuve}
$R_q(1,m)$ is self-complementary, so, by theorem \ref{gen}, it is enough to show that $R_q(1,m)$ has strength 2.
\\Let $y=(y_1,\ldots,y_m)$ and $z=(z_1,\ldots,z_m)$ two different fixed elements of $\F{q}^m$.
\\Let $f$, $g\in R_q(1,m)$, we say that $f$ is equivalent to $g$ ($f\sim g$) if and only if $f(y)=g(y)$ and $f(z)=g(z)$.
\\Let $f(x)=a_1x_1+\ldots+a_mx_m+b$. Let $g=\alpha_1x_1+\ldots+\alpha_mx_m+\beta$ such that $f \sim g$, then
\begin{eqnarray*}&&\left\{\begin{array}{l}\alpha_1y_1+\ldots+\alpha_my_m+\beta=a_1y_1+\ldots+a_my_m+b\\
\alpha_1z_1+\ldots+\alpha_mz_m+\beta=a_1z_1+\ldots+a_mz_m+b\end{array}\right.\\&&\qquad\qquad\Leftrightarrow\left(\begin{array}{cccc}y_1& \ldots&y_m&1\\z_1&\ldots&z_m&1\end{array}\right)\left(\begin{array}{c}\alpha_1\\ \vdots\\ \alpha_m\\ \beta\end{array}\right)=\left(\begin{array}{c}a_1y_1+\ldots+a_my_m+b\\a_1z_1+\ldots+a_mz_m+b\end{array}\right)\end{eqnarray*}
Since $z\neq y$, this system has rank 2, and so there are $q^{m-1}$ solutions.
\\Furthermore, $\mathrm{Card}(R_q(1,m))=q^{m+1}$, and $\frac{q^{m+1}}{q^{m-1}}=q^2=\mathrm{Card}(\mathbb{F}_q^2)$.
\\Hence $R_q(1,m)$ has strength 2.
\end{preuve}

\begin{proposition}\label{min}For all $q$, $m\geq 0$, we have $$\rho(1,m)\geq (q-1)q^{m-1}-q^{\lceil \frac{m}{2}\rceil-1}.$$\end{proposition}

\begin{preuve} We have $$\rho(1,m)\geq\rho_2(1,m),$$
and by theorem \ref{rho2}, $$\rho_2(1,m)=(q-1)q^{m-1}-q^{\lceil \frac{m}{2}\rceil-1},$$
which gives the result.\end{preuve}

\begin{remarque}$\rho(1,1)=q-2$.
\\Indeed, by proposition \ref{min}, $\rho(1,1)\geq q-2$.
Furthermore, for all $g\in B_1^q$, we consider $f(x)=g(x)-(ax+b)$ with $a=g(1)-g(0)$ and $b=g(0)$; $f$ has at least two roots (0 and 1) so for all $g\in B_1^q$, $d(g,R_q(1,1))\leq q-2$. Hence $\rho(1,1)\leq q-2$.\end{remarque}

Combining proposition \ref{max} and \ref{min} we get the following :

\begin{corollaire}For all $q$, if $m$ is even, then
$$\rho(1,m)=(q-1)q^{m-1}-q^{\frac{m}{2}-1}.$$\end{corollaire}

\begin{preuve}For $m$ even, $(q-1)q^{m-1}-q^{\lceil\frac{m}{2}\rceil-1}=(q-1)q^{m-1}-q^{\frac{m}{2}-1}$.\end{preuve}

\begin{remarque} Furthermore, we have shown that $\rho_2(1,m)=(q-1)q^{m-1}-q^{\lceil \frac{m}{2}\rceil-1}$ and that $\rho_2(1,m)$  is reached for $f(x)=x_1x_2+\ldots+x_{m-3}x_{m-2}+ax_{m-1}^2+bx_m^2+cx_{m-1}x_m$ with $ax^2+cx+b$ irreducible over $\mathbb{F}_q$.
 So, since for $m$ even $\rho(1,m)=\rho_2(1,m)$, we get that $\rho(1,m)=d(f,R_q(1,m))$.\end{remarque}

\section{Calculation of $\rho(1,3)$ for $q=3$}

From now, we assume that $q=3$.

\begin{theoreme}\label{rho3} For $q=3$, $\rho(1,3)=16$.\end{theoreme}

\begin{preuve}By theorem \ref{bound}, $15\leq\rho(1,3)\leq16$. Furthermore, if $\rho(1,3)=16$, there exists $f\in B_m^q$ such that $d(f,R_3(1,3))=16$ and necessarily, degree of $f$ is greater than 2 since, by theorem \ref{rho2}, $\rho_2(1,3)=15$.
\\Using all these restrictions, we use Magma
\begin{algo}{}\begin{verbatim}K:=GF(3);
P<x,y,z>:=PolynomialRing(K,3);
R1:=[a*x+b*y+c*z+d : a in K, b in K, c in K, d in K];
M:=15;
L:=[0,0,0,0,0,0,1,0,0,0,0,0,0,0,0,0,0,0,0,0,0,0,0];

ad:=function(L); i:=1; r:=true;
while i le #L and r do
if L[i]+1 eq 3 then i:=i+1; L[i-1]:=0;
   else r:=false; L[i]:=L[i]+1; end if; end while;
   return L;
 end function;

while M eq 15  and L[23] le 1 do
pol:=L[1]*z^2+L[2]*y*z+L[3]*y^2+L[4]*x*z+L[5]*x*y+L[6]*x^2
 +L[7]*y*z^2+L[8]*y^2*z+L[9]*x*z^2+L[10]*x*y*z+L[11]*x*y^2
 +L[12]*x^2*z+L[13]*x^2*y+L[14]*y^2*z^2+L[15]*x*y*z^2
 +L[16]*x*y^2*z+L[17]*x^2*z^2+L[18]*x^2*y*z+L[19]*x^2*y^2
 +L[20]*x*y^2*z^2+L[21]*x^2*y*z^2+L[22]*x^2*y^2*z
 +L[23]*x^2*y^2*z^2;
k:=1; m:=16;
while k le #R1 and m eq 16 do
if Evaluate(R1[k],<0,0,0>) ne 0 then r:=1; else; r:=0; end if;
p:=<1,0,0>;
while r lt 16 and p ne <0,0,0>  do
if Evaluate(pol+R1[k],p) ne 0 then r:=r+1; p:=ad(p);
   else p:=ad(p); end if; end while;
if r lt m then m:=r; end if;
k:=k+1; end while;
if m gt M then M:=m; else L:=ad(L); end if; end while;

print(M);
print(L);
\end{verbatim}\end{algo}
We get that $d(y^2+xy+y^2z+xyz+y^2z^2+x^2z^2,R_3(1,3))=16$.\end{preuve}

\begin{proposition}There is no $f\in R_3(6,4)\setminus R_3(4,3)$ such that $d(f,R_3(1,3))=16$.\end{proposition}

\begin{lemme}\label{deg6}For $q\geq3$, if $f\in R_q(m(q-1),m)\setminus R_q(m(q-1)-1,m)$ then there exists $\sigma\in \mathrm{GA}_m(\F{q})$, $a\in\F{q}^*$ and $r\in R_q(m(q-1)-2,m)$ such that $$\sigma .f=a\prod_{i=1}^m x_i^{q-1} +r$$\end{lemme}

\begin{preuve}
We write $f$ as $$f=a\prod_{i=1}^m x_i^{q-1}+\sum_{i=1}^m a_i \prod_{k\neq i}x_k^{q-1}x_i^{q-2}+s$$ where $a$, $a_i\in\F{q}$, $a\neq0$ and $s\in R_q(m(q-1)-2,m)$.
\\\\Let $\omega\in\F{q}^m$ then 
\begin{eqnarray*}\mathrm{1}_\omega&=&\prod_{i=1}^m(1-(x_i-\omega_i)^{q-1}) \\&=&\prod_{i=1}^m\left(1-\sum_{k=1}^{q-1}\binom{q-1}{k}x_i^k(-\omega_i)^{q-1-k}\right) \\&=&(-1)^m\prod_{i=1}^mx_i^{q-1}+(-1)^m\sum_{i=1}^m\omega_i\prod_{k\neq i}x_k^{q-1}x_i^{q-2}+t,\qquad\textrm{$t\in R_q(m(q-1)-2,m)$}.\end{eqnarray*} 
Hence $$f=(-1)^ma\mathrm{1}_{(a^{-1}a_i)}+r', \qquad r'\in R_q(m(q-1)-2,m).$$
Let $\sigma\in\mathrm{GA}_m(\F{q})$, $$\sigma.f=(-1)^ma\mathrm{1}_{\sigma^{-1}(a^{-1}a_i)}+\sigma.r',\qquad r'\in R_q(m(q-1)-2,m).$$
We choose $\sigma$ such that $\sigma^{-1}(a^{-1}a_i)=0$.
\begin{eqnarray*}\mathrm{1}_0&=&\prod_{i=1}^m(1-x_i^{q-1})=(-1)^m\prod_{i=1}^mx_i^{q-1}+u,\\&&\hspace{4cm} u\in R_q((m-1)(q-1),m)\subset R_q(m(q-1)-2,m)\quad \textrm{since $q\geq 3$}\end{eqnarray*}
Finally, since $\mathrm{Aut}(R_q(m(q-1)-2,m))=\mathrm{GA}_m(\F{q})$, we get $$\sigma.f=a\prod_{i=1}^mx_i^{q-1}+r,\qquad r\in R_q(m(q-1)-2,m).$$
\end{preuve}

\begin{lemme}\label{deg5}If $f\in R_q(m(q-1)-1,m)\setminus R_q(m(q-1)-2,m)$ then there exists $\sigma\in \mathrm{GL}_m(\F{q})$ and $r\in R_q(m(q-1)-2,m)$ such that $$\sigma .f=\prod_{i=1}^{m-1} x_i^{q-1}x_m^{q-2} +r$$\end{lemme}

\begin{preuve}We write $f=\displaystyle\sum_{i=1}^m\alpha_i\prod_{k\neq i}x_k^{q-1}x_i^{q-2}+t$, $t\in R_q(m(q-1)-2,m)$.
\\\\Let $b \in \F{q}^m$,
\begin{eqnarray*}\mathrm{1}_0 -\mathrm{1}_b&=&\prod_{i=1}^m(1-x_i^{q-1})-\prod_{i=1}^m(1-(x_i-b_i)^{q-1})\\&=&\prod_{i=1}^m(1-x_i^{q-1})-\prod_{i=1}^m(1-\sum_{k=1}^{q-1}\binom{q-1}{k}x_i^k(-b_i)^{q-1-k})\\&=&(-1)^{m+1}\sum_{i=1}^mb_i\prod_{k\neq i}x_k^{q-1}x_i^{q-2}+s,\qquad\textrm{$s\in R_q(m(q-1)-2,m)$.} \end{eqnarray*}
D'où $f=\mathrm{1}_0 -\mathrm{1}_{((-1)^{m+1}\alpha_i)}+r'$, $r'\in R_q(m(q-1)-2,m).$
\\\\Let $\sigma\in\mathrm{GL}_m(\F{q})$ then 
$$\sigma.f=\mathrm{1}_{\sigma^{-1}(0)} -\mathrm{1}_{\sigma^{-1}((-1)^{m+1}\alpha_i)}+\sigma.r'=\mathrm{1}_0 -\mathrm{1}_{\sigma^{-1}((-1)^{m+1}\alpha_i)}+\sigma.r'$$
Since $f\in R_q(m(q-1)-1,m)\setminus R_q(m(q-1)-2,m)$, $((-1)^{m+1}\alpha_i)\neq 0$. So there exists $\sigma \in \mathrm{GL}_m(\F{q})$ such that $$\sigma^{-1}((-1)^{m+1}\alpha_i)=\left(\begin{array}{c}(-1)^{m+1}\\ 0 \\ \vdots \\ 0\end{array}\right)=c$$
and $$\sigma.f=\mathrm{1}_0 -\mathrm{1}_c+\sigma.r'=\displaystyle\prod_{i=1}^{m-1}x_i^{q-1}x_m^{q-1}+r$$
with $r\in R_q(m(q-1)-2,m)$ since $\mathrm{Aut}(R_q(m(q-1)-2,m)=\mathrm{GA}_m(\F{q})$.\end{preuve}

\begin{preuvep}By lemma \ref{deg6} and \ref{deg5}, the following algorithms on Magma give the result.
\begin{algo}{(degré 6)}\label{degre6}
\begin{verbatim}K:=GF(3);
P<x,y,z>:=PolynomialRing(K,3);
R1:=[a*x+b*y+c*z+d : a in K, b in K, c in K, d in K];

ad:=function(L); i:=1; r:=true; 
while i le #L and r do 
   if L[i]+1 eq 3 then i:=i+1; L[i-1]:=0; 
      else r:=false; L[i]:=L[i]+1; end if; end while; 
return L;
end function;

L:=[1,0,0,0,0,0,0,0,0,0,0,0,0,0,0,0,0,0,0];
while L ne [0,0,0,0,0,0,0,0,0,0,0,0,0,0,0,0,0,0,0] do  
pol:=L[1]*z^2+L[2]*y*z+L[3]*y^2+L[4]*x*z+L[5]*x*y+L[6]*x^2
     +L[7]*y*z^2+L[8]*y^2*z+L[9]*x*z^2+L[10]*x*y*z+L[11]*x*y^2
     +L[12]*x^2*z+L[13]*x^2*y+L[14]*y^2*z^2+L[15]*x*y*z^2
     +L[16]*x*y^2*z+L[17]*x^2*z^2+L[18]*x^2*y*z
     +L[19]*x^2*y^2+x^2*y^2*z^2;
k:=1; m:=16;
while k le #R1 and m eq 16 do
if Evaluate(R1[k],<0,0,0>) ne 0 then r:=1; else; r:=0; end if;
p:=<1,0,0>;
while r lt 16 and p ne <0,0,0>  do 
  if Evaluate(pol+R1[k],p) ne 0 then r:=r+1; p:=ad(p); 
   else p:=ad(p); end if; end while;
if r lt m then m:=r; end if;
k:=k+1; end while; 
if m gt 15 then print(pol); L:=ad(L);  else L:=ad(L); end if;  end while;
\end{verbatim}\end{algo}
\newpage
\begin{algo}{(degré 5)}\label{degre5}\begin{verbatim}K:=GF(3);
P<x,y,z>:=PolynomialRing(K,3);
R1:=[a*x+b*y+c*z+d : a in K, b in K, c in K, d in K];

ad:=function(L); i:=1; r:=true; 
while i le #L and r do 
 if L[i]+1 eq 3 then i:=i+1; L[i-1]:=0; 
   else r:=false; L[i]:=L[i]+1; end if; end while; 
 return L;
 end function;
 
L:=[1,0,0,0,0,0,0,0,0,0,0,0,0,0,0,0,0,0,0];
while L ne [0,0,0,0,0,0,0,0,0,0,0,0,0,0,0,0,0,0,0] do  
pol:=L[1]*z^2+L[2]*y*z+L[3]*y^2+L[4]*x*z+L[5]*x*y+L[6]*x^2
     +L[7]*y*z^2+L[8]*y^2*z+L[9]*x*z^2+L[10]*x*y*z+L[11]*x*y^2
     +L[12]*x^2*z+L[13]*x^2*y+L[14]*y^2*z^2+L[15]*x*y*z^2
     +L[16]*x*y^2*z+L[17]*x^2*z^2+L[18]*x^2*y*z+L[19]*x^2*y^2+x^2*y^2*z;
k:=1; m:=16;
while k le #R1 and m eq 16 do
if Evaluate(R1[k],<0,0,0>) ne 0 then r:=1; else; r:=0; end if;
p:=<1,0,0>;
while r lt 16 and p ne <0,0,0>  do 
 if Evaluate(pol+R1[k],p) ne 0 then r:=r+1; p:=ad(p); 
   else p:=ad(p); end if; end while;
if r lt m then m:=r; end if;
k:=k+1; end while;
if m gt 15 then print(pol); L:=ad(L);  else L:=ad(L); end if; end while;\end{verbatim}\end{algo}
Both algorithms do not give any $f$ such that $d(f,R_3(1,3))=16$.
\end{preuvep}

Using a similar algorithm, we can verify the following proposition :

\begin{proposition}All $f$ in $R_3(4,3)$ such that $d(f,R_3(1,3))=16$ are equivalent under the action of $\mathrm{GA}_3(\F{3})$ and $R_3(1,3)$ to $$2x^2z^2+2yz+x^2z^2+xyz+2x^2yz.$$\end{proposition}

\section{Improvement of the lower bound of $\rho(1,m)$ for $q=3$}

We use theorem \ref{rho3} to improve, for $q=3$,  the lower bound of $\rho(1,m)$ given by theorem \ref{bound}.

\begin{lemme}\label{rec}For all $q$, for all $m$,
$$\rho(1,m+2)\geq (q-1)^2q^m+q\rho(1,m)$$\end{lemme}

\begin{preuve} Let $f\in R_q(1,m+2)$. We can write $f(x_1,\ldots,x_{m+2})=g(x_1,\ldots,x_m)+\alpha x_{m+1}+\beta x_{m+2}$, where $g\in R_q(1,m)$ and $\alpha$, $\beta\in\F{q}$. \\We denote the elements of $\F{q}$ by $\omega_0$, $\omega_1,\ldots, \omega_{q-1}$. We can assume that $\omega_0=0$.
Then $$R_q(1,m+2)=\{M(c,\alpha,\beta), c\in R_q(1,m),\alpha\in\mathbb{F}_q,\beta\in \mathbb{F}_q\}$$
where  {$M(c,\alpha,\beta)=(c,\overline{c}^{\alpha\omega_1},\ldots,\overline{c}^{\alpha\omega_{q-1}},\overline{c}^{\beta\omega_1}, \overline{c}^{\alpha\omega_1+\beta\omega_1},\ldots,\overline{c}^{\beta\omega_1+\alpha\omega_{q-1}},\ldots,\overline{c}^{\alpha\omega_{q-1}+ \beta\omega_{q-1}})$}
\\\\If $v_0$ is such that $d(v_0,R_q(1,m))=\rho(1,m)$, then for all $c\in R_q(1,m)$, $|v_0+c|\geq\rho(1,m)$.
\\Let $u=(\overline{v_0}^{\omega_0\omega_0},\ldots,\overline{v_0}^{\omega_0\omega_{q-1}},\overline{v_0}^{\omega_1\omega_{0}},\ldots,\overline{v_0}^{\omega_1\omega_{q-1}}, \ldots,\overline{v_0}^{\omega_{q-1}\omega_{q-1}})\in B_{m+2}^q$. \\If $\alpha$, $\beta\in\mathbb{F}_q$ and $c\in R_q(1,m)$, then
\begin{eqnarray*}|u+M(c,\alpha,\beta)|&=&\sum_{i=0}^{q-1}\sum_{j=0}^{q-1}|\overline{v_0}^{\omega_i\omega_j}+\overline{c}^{\alpha\omega_i+\beta\omega_j}|
\\&=&\sum_{i=0}^{q-1}\sum_{j=0}^{q-1}|\overline{v_0}^{(\beta+\omega_i)\omega_j}+\overline{c}^{\alpha\omega_i}|\\&=&q|\overline{c}^{\alpha(-\beta)}+v_0|+ \sum_{\omega_i\neq-\beta}\sum_{j=0}^{q-1}|\overline{v_0}^{(\beta+\omega_i)\omega_j}+\overline{c}^{\alpha\omega_i}|\\&\geq&q\rho(1,m)+(q-1)^2q^m\end{eqnarray*}
which gives the result.
\end{preuve}

\begin{theoreme}\label{min2}For $q=3$ and $m$ an odd integer, $$\rho(1,m)\geq (q-1)q^{m-1}-\frac{2}{3}q^{\lceil\frac{m}{2}\rceil-1}$$\end{theoreme}

\begin{preuve} We write $m=2k+1$. We prove by induction on $u$ that for $u\leq k$, $$\rho(1,m)\geq(q-1)(q^{m-1}-q^{m-1-u})+q^u\rho(1,m-2u).$$
This is true for $u=0$ and $u=1$ (lemma \ref{rec}). Assume that it is true for some $u<k$. Then
{\setlength\arraycolsep{2pt}\begin{eqnarray*}\rho(1,m)&\geq&(q-1)(q^{m-1}-q^{m-1-u})+q^u\rho(1,m-2u) \\&\geq&(q-1)(q^{m-1}-q^{m-1-u})+q^u\left((q-1)^2q^{m-2u-2}+q\rho(1,m-2u-2)\right) \\&&\hspace{8cm} \textrm{(by lemma \ref{rec})} \\&\geq&(q-1)q^{m-1}-(q-1)(q-(q-1))q^{m-u-2}+q^{u+1}\rho(1,m-2(u+1)) \\&\geq&(q-1)(q^{m-1}-q^{m-(u+1)-1})+q^{u+1}\rho(1,m-2(u+1))\end{eqnarray*}}
Hence, for $q=3$ and $ u=k-1$, we get :
{\setlength\arraycolsep{2pt}\begin{eqnarray*}\rho(1,m)&\geq&(q-1)(q^{m-1}-q^{k+1})+q^{k-1}\rho(1,3) \\&\geq&(q-1)q^{m-1}-2q^{k-1} =(q-1)q^{m-1}-\frac{2}{3}q^{\lceil\frac{m}{2}\rceil-1}\end{eqnarray*}}
\end{preuve}

\begin{corollaire}For $q=3$, $\rho(1,5)=156$\end{corollaire}

\begin{preuve}By theorem \ref{min2}, $\rho(1,5)\geq 2.3^4-\frac{2}{3}.9=156$.
and by proposition \ref{max}, $\rho(1,5)\leq[2.3^4-3\sqrt{3}]=156$.\end{preuve}

\end{document}